\input amstex
\documentstyle{amsppt}
\NoBlackBoxes
\mag1100
\document
\topmatter
\title 
Correction to  ``A smooth foliation of the 5-sphere by complex surfaces''.
\endtitle 
\rightheadtext{corrigendum}
\author
Laurent Meersseman,
Alberto Verjovsky
\endauthor
\date June 2, 2011\enddate 
\address 
{Laurent Meersseman}\hfill\hfill\linebreak
\indent{I.M.B.}\hfill\hfill\linebreak
\indent{Universit\'e de Bourgogne}\hfill\hfill\linebreak
\indent{B.P. 47870}\hfill\hfill\linebreak
\indent{21078 Dijon Cedex}\hfill\hfill\linebreak
\indent{France}\hfill\hfill
\endaddress
\email laurent.meersseman\@u-bourgogne.fr \endemail

\address
{Alberto Verjovsky}\hfill\hfill\linebreak \indent{Instituto de
Matem\'aticas de la UNAM, Unidad Cuernavaca}\hfill\hfill\linebreak
\indent{Apartado Postal 273-3, Admon. de correos
No.3}\hfill\hfill\linebreak \indent{Cuernavaca, Morelos,
M\'exico}\hfill\hfill
\endaddress
\email alberto\@matcuer.unam.mx \endemail
\endtopmatter

Our paper \cite{M-V} claims that it describes a foliation of $\Bbb S^5$ by complex surfaces. However it was pointed out to us by the anonymous referee of a related article that the foliation constructed in the paper lives in fact on a $5$-manifold with non-trivial fundamental group. The aim of this note is to explain this difference and to characterize this $5$-manifold.
\medskip
We observe that, even with this modification, this foliation is still the first example of such an exotic CR-Structure. Quoting \cite{M-V2}, ``We would like to emphasize that, as far as we know, the foliation described
in [M-V] (as well as the related examples of [M-V], Section 5) is the only known example of
a smooth foliation by complex manifolds of complex dimension strictly greater than one on
a compact manifold, which is not obtained by classical methods such as the one given by
the orbits of a locally free smooth action of a complex Lie group, the natural product foliation
on $M\times N$ where $M$ is foliated by Riemann surfaces and $N$ is a complex manifold,
holomorphic fibrations, or trivial modifications of these examples such as cartesian products
of known examples or pull-backs. Of course, it is very easy to give examples of foliations
by complex manifolds on open manifolds (in fact even with Stein leaves). On the other
hand, if a compact smooth manifold has an orientable smooth foliation by surfaces then,
using a Riemannian metric and the existence of isothermal coordinates, we see that the foliation
can be considered as a foliation by Riemann surfaces."
\medskip
We use the notations and results of \cite{M-V} and assume that the reader is acquainted with them.
\medskip
The foliation of \cite{M-V} is obtained by gluing, thanks to Lemma 1, two tame foliations on manifolds with boundary. The first one, $\Cal M$ is a bundle over the circle with fiber the affine Fermat surface
$$
F=\{(z_1,z_2,z_3)\in\Bbb C^3\quad\vert\quad P(z)=z_1^3+z_2^3+z_3^3=1\}
$$
and monodromy
$$
z\in F\longmapsto \omega\cdot z\in F
$$
where $\omega=\exp (2i\pi/3)$. Its foliation is described in \cite{M-V, Section 3}.
\medskip
The second one, $\Cal N$, is supposed to be diffeomorphic to $K\times \overline{\Bbb D}$ (where $K$ is  a circle bundle over a torus, see \cite{M-V, Section 1.2}). Its foliation is given as a quotient foliation: take the quotient of 
$$
\tilde X=\Bbb C^*\times (\Bbb C\times [0,\infty )\setminus\{(0,0)\}),
$$
endowed with the trivial foliation by the level sets of the $[0,\infty )$-coordinate, by the abelian group generated by
$$
\left\{
\eqalign{T(z,u,t)&=(\exp (2i\pi \omega)\cdot z, ((\psi(z))^{-3}\cdot u, t)\cr
U(z,u,t)&=(z,\exp (2i\pi\tau) \cdot u, d(t))
}
\right .
$$
(cf. \cite{M-V, Section 2} for the definition of $\psi$ and $d$). The map $T$ is chosen such that the quotient space of the boundary of $\tilde X$ (corresponding to $t=0$) by the $T$-action is a particular $\Bbb C^*$ bundle of Chern class $-3$ over the elliptic curve $\Bbb E_\omega$, called $W$ in \cite{M-V}.
\medskip
Consider the interior of this manifold with boundary (that is take $t>0$). The quotient of a leaf $\{t=Constant\}$ by the $T$-action is biholomorphic to $L$, the line bundle over $\Bbb E_\omega$ associated to $W$. Hence the quotient of $\text{Int }\tilde X$ by the $T$-action is CR-isomorphic to $L\times (0,\infty )$. Now, using the fact that $d$ is contracting, and the fact that the map
$$
u\in\Bbb C\longmapsto \exp (2i\pi\tau) \cdot u\in\Bbb C
$$
is isotopic to the identity, we see that the complete quotient is a bundle over the circle with fiber $L$ and monodromy isotopic to the identity, that is CR-isomorphic to $L\times\Bbb S^1$.
Indeed, this means that $\Cal N$ (if we decide from now on to call $\Cal N$ the previous quotient manifold with boundary) is diffeomorphic to $D\times \Bbb S^1$, where $D$ is the $\overline{\Bbb D}$-bundle associated to $L$. It is definitely not diffeomorphic to $K\times 
  \overline{\Bbb D}$, since this last manifold has a nilpotent fundamental group (see \cite{M-V, Section 1.2}), whereas $D\times \Bbb S^1$ retracts on $\Bbb S^1\times\Bbb S^1\times\Bbb S^1$.
\medskip
When gluing $\Cal M$ and this ``new" $\Cal N$, one does not obtain the $5$-sphere but the following manifold, let us call it $Z$. Observe that, in the projective Fermat surface
$$
F^p=\{[z_0:z_1:z_2:z_3]\in\Bbb P^3\quad\vert\quad z_1^3+z_2^3+z_3^3=z_0^3\},
$$
the elliptic curve at infinity has a tubular neighborhood diffeomorphic to $D$. Observe also that the gluing of $\Cal M$ and $\Cal N$ respect the fibrations over the circle, that is the two bases are identified and the gluing occur on the fibers. It follows from all that that $Z$ is a bundle over the circle with fiber $F^p$ and monodromy
$$
[z_0:z_1:z_2:z_3]\in F^p\longmapsto [z_0:\omega\cdot z_1:\omega\cdot z_2:\omega\cdot z_3]\in F^p
$$

So finally, what is really proved in \cite{M-V} is the following Theorem.

\proclaim{Theorem}
Let $Z$ be the $5$-dimensional bundle over the circle with fiber the projective Fermat surface and monodromy the multiplication by the root of unity $\omega$ on the affine part.

There exists on $Z$ an exotic smooth, codimension-one, integrable and Levi-flat CR-structure on $Z$. The induced foliation by complex surfaces satisfies:

(i) There are only two compact leaves both biholomorphic to an elliptic bundle over the elliptic curve $\Bbb E_\omega$. Since this surface has odd first Betti number it is not K\"ahler.

(ii) One compact leaf is the boundary of a compact set in $Z$ whose interior is foliated by line bundles over $\Bbb E_\omega$ with Chern class $-3$ obtained from
$W$ by adding a zero section. The two compact leaves are the boundary
components of a collar and the leaves in the interior of this collar are
infinite cyclic coverings of the compact leaves and biholomorphic to $W$,
thus they are principal $\Bbb C^*$-bundles over the elliptic curve $\Bbb E_\omega$.

(iii) The other leaves have the homotopy type of a bouquet of eight copies of
$\Bbb S^2$ and they are all biholomorphic to the affine complex smooth manifold
$P^{-1} ({z}),\ z \in\Bbb C^*$.
\endproclaim

\remark{Remarks}

\noindent a) Since the manifold $Z$ fibres over the circle with fibre a
nonsigular cubic surface, it has a natural foliation by complex
leaves
given by the fibres. Ours is obviously completely different.

\noindent b) Using the polynomials
$$
P(z)=z_1^2+z_2^4+z_3^4\qquad\text{resp. }\qquad P(z)=z_1^2+z_2^3+z_3^6
$$
instead of the cubic one, the construction can be adapted as described in \cite{M-V} to obtain exotic integrable CR-structures on bundles over the circle with fiber
$$
F^p=\{[z_0:z_1:z_2:z_3]\in\Bbb P^3\quad\vert\quad z_0^2z_1^2+z_2^4+z_3^4=z_0^4\},
$$
and respectively
$$
F^p=\{[z_0:z_1:z_2:z_3]\in\Bbb P^3\quad\vert\quad z_0^4z_1^2+z_0^3z_2^3+z_3^6=z_0^6\}.
$$

\noindent c) Due to the compact non K\"ahler leaves, this CR-structure is not embeddable in any Stein space nor K\"ahler manifold. Moreover, it is not embeddable in any $3$-dimensional complex manifold \cite{DS}.

\noindent d) In \cite{De}, G. Deschamps proved that the use of a collar can be avoided by choosing carefully the holonomy of the boundaries of $\Cal M$ and $\Cal N$. This gives a foliation on $Z$ with the same properties as above, but with a single compact leaf.
\endremark


\vskip1cm
\Refs
\widestnumber\key{ooooo}
\ref
\key De
\by G. Deschamps
\paper Feuilletage lisse de $\Bbb S^5$ par surfaces complexes
\jour C.R. Acad. Sci. Paris \vol 348 \yr 2010 \pages 1303--1306
\endref

\ref
\key DS
\by G. Della Sala
\paper Non-embeddability of certain classes of Levi flat manifolds
\paperinfo preprint
\yr 2009
\endref

\ref
\key M-V
\by L. Meersseman, A. Verjovsky
\paper A smooth foliation of the 5-sphere by complex surfaces
\jour Ann. Math. \vol 156 \yr 2002 \pages 915--930
\endref

\ref
\key M-V2
\by L. Meersseman, A. Verjovsky
\paper On the moduli space of certain smooth codimension-one foliations of the $5$-sphere
\jour J. Reine Angew. Math.
\vol 632 \pages 143--202
\yr 2009
\endref
\endRefs
\enddocument